\documentclass{amsart}
\usepackage{amsmath, amssymb,verbatim}
\newtheorem{theorem}{Theorem}
\newtheorem{proposition}{Proposition}

\newtheorem{lemma}{Lemma}
\newtheorem{remark}{Remark}

\subjclass[2010]{Primary 47B35; Secondary 47L80}
\keywords{Toeplitz operators, quasihomogeneous symbol, Mellin transform.}
\author[BOUHALI]{Aissa Bouhali}
\address{D\'epartement de Math\'ematiques, \'Ecole Normale Sup\'erieure de Laghouat;
	Laboratory of Pure and Applied Mathematics, University of Laghouat; Applied Sciences and Didactics Laboratory, ENS of Laghouat, Algeria.}
\email{aissa.bouhali@ens-lagh.dz}
\author[LOUHICHI]{Issam Louhichi }
\address{Department of Mathematics \& Statistics, College of Arts \& Sciences, American University of Sharjah, P.O.Box 2666, Sharjah, UAE.}
\email{ilouhichi@aus.edu}
\author[Yousef]{Abdel Rahman Yousef }
\address{Department of Mathematics \& Statistics, College of Arts \& Sciences, American University of Sharjah, P.O.Box 2666, Sharjah, UAE.}
\email{afyousef@aus.edu}

\begin{document}
	\title[Commutants of a certain class of Toeplitz operators ]
	{Commutants of a certain class of  Toeplitz operators}

	\date{\today} 
	\begin{abstract} A major open problem in the Theory of Toeplitz operators on the analytic Bergman space over the unit disk is the characterization of the commutant of a given Toeplitz operator--that is, the set of all bounded Toeplitz operators that commute with it. In this paper, we provide a complete description of bounded Toeplitz operators $T_f$, where the symbol $f$ has a truncated polar decomposition, that commute with a Toeplitz operator whose symbol is the sum of a quasihomogeneous function and a bounded analytic function.
\end{abstract}
\maketitle 
\section{Introduction}
In the complex plane \(\mathbb{C}\), let \(\mathbb{D}\) denote the open unit disk, and consider the normalized Lebesgue measure \(dA = r \, dr \, \frac{d\theta}{\pi}\) on \(\mathbb{D}\), where \((r, \theta)\) are the polar coordinates. The space \(L^2(\mathbb{D}, dA)\) consists of all square-integrable functions on \(\mathbb{D}\) with respect to this measure.  

The classical unweighted Bergman space, denoted by \(L^2_a(\mathbb{D})\), is the closed subspace of \(L^2(\mathbb{D}, dA)\) comprising all functions that are analytic in \(\mathbb{D}\). The set \(\{z^k : k = 0, 1, 2, \dots\}\) forms an orthogonal basis for \(L^2_a(\mathbb{D})\). Since this space is closed, there exists a well-defined orthogonal projection \(P\) from \(L^2(\mathbb{D}, dA)\) onto \(L^2_a(\mathbb{D})\), commonly referred to as the Bergman projection. For a comprehensive treatment of Bergman spaces and their associated projections, readers may refer to \cite{h}.  

Given a function \(f \in L^1(\mathbb{D}, dA)\), the Toeplitz operator \(T_f\) acting on \(L^2_a(\mathbb{D})\) is defined by \(T_f(g) = P(fg)\), provided that \(fg \in L^2(\mathbb{D}, dA)\). Here, \(f\) is referred to as the symbol of the Toeplitz operator \(T_f\). This definition implies that any bounded analytic function on \(\mathbb{D}\) belongs to the domain of \(T_f\), ensuring that \(T_f\) is densely defined. Moreover, if the symbol \(f\) is bounded on \(\mathbb{D}\), the Toeplitz operator \(T_f\) is also bounded, with the norm satisfying \(\|T_f\| \leq \|f\|_\infty\). Nevertheless, an unbounded symbol $f$ belonging to $L^1(\mathbb{D},DA)$ may still generate a bounded Toeplitz operator $T_f$. In fact, if $f$ is integrable and remains bounded on the annulus $\{z\in\mathbb{D}\ :\ 0<r<|z|<1\}$ for some $0<r<1$, then $f$ can be decomposed into the sum of two functions: one that is integrable with compact support and another that is bounded. Consequently, $T_f$ remains bounded. Such symbols are referred to as "nearly bounded symbols" \cite[~p.204]{ahern}. From now on, we will consider any symbol that  does not belong $L^1(\mathbb{D},dA)$ as inadmissible for a bounded Toeplitz operator, meaning such a symbol cannot define a bounded Toeplitz operator.

A function \(f\) is called quasihomogeneous of degree \(p\), where \(p\) is an integer, if it can be expressed in the form \(f(re^{i\theta}) = e^{ip\theta}\phi(r)\), where \(\phi\) is a radial function. The Toeplitz operator \(T_f\) associated with such a symbol is then referred to as a quasihomogeneous Toeplitz operator of degree \(p\) (see \cite{cr, lsz}). The study of these operators is motivated by the structural decomposition of \(L^2(\mathbb{D}, dA)\), which can be written as \(L^2(\mathbb{D}, dA) = \bigoplus_{k \in \mathbb{Z}} e^{ik\theta} \mathcal{R}\), where \(\mathcal{R}\) denotes the space of square-integrable radial functions on \([0,1)\) with respect to the measure \(r \, dr\).  This decomposition implies that any function \(f \in L^2(\mathbb{D}, dA)\) admits a polar expansion \(f(z) = f(re^{i\theta}) = \sum_{k \in \mathbb{Z}} e^{ik\theta} f_k(r)\), where each \(f_k(r)\) is a radial function. The subclass of quasihomogeneous Toeplitz operators provides a natural framework for understanding the behavior of operators that respect this decomposition (see \cite{bl, cr, l, lr, lry, lsz, lz, sz}).

Studying this family of operators offers insight into the commutant of Toeplitz operators associated with more general symbols. By focusing on quasihomogeneous symbols, we aim to uncover structural properties and extend our understanding of operator commutativity in the Bergman space setting.   

Over the past several decades, Toeplitz operators have been a central object of study, particularly in terms of their algebraic and functional properties. One of the most intriguing open problems involves characterizing the commutant of a Toeplitz operator, i.e., determining the set of all Toeplitz operators that commute with a given one under composition. Despite substantial progress, little is known about the commutativity of Toeplitz operators for general symbols \(f\) and \(g\). For an overview of known results on the commutativity of Toeplitz operators in \(L^2_a(\mathbb{D})\) and specific cases where progress has been made, see \cite{sl, al, acr,bl, cr, lt, l, lr, lry, lsz, lz, rv, sz, y}.  
	
	\section{Tools}
	The Mellin transform $\widehat{\phi}$ of a radial function $\phi\in L^{1}([0,1), rdr)$ is defined as $$\widehat{\phi}(z)=\int_0^1\phi(r)r^{z-1}dr.$$
	It is well-known that for such functions, the Mellin transform is bounded  on the right-halfplane $\{z:\Re z\geq 2\}$ and is analytic on $\{z:\Re z> 2\}$. 
	
		The following lemma gives the action of quasihomogeneous Toeplitz operators on the vectors of the orthogonal basis of $L^2_a(\mathbb{D})$. 
	\begin{lemma}\label{mellin}\cite[Lemma 5.3, p.531]{lsz}
		Let $k, p\in\mathbb{Z}_+$ and let $\phi$ be a radial function in $L^1([0,1), rdr)$. Then
		$$T_{e^{ip\theta}\phi}(z^k)=2(k+p+1)\widehat{\phi}(2k+p+2)z^{k+p}$$ and 
		$$T_{e^{-ip\theta}\phi}(z^k)=\left\{\begin{array}{ll}0&\textrm{ if }\ 0\leq k\leq p-1,\\
			2(k-p+1)\widehat{\phi}(2k-p+2)z^{k-p}&\textrm{ if }\ k\geq p.\end{array}\right.$$
	\end{lemma}	
	The Mellin transform of a function is uniquely determined by its values on any arithmetic sequence of integers.  
	In fact we have the following classical theorem \cite[p.102]{Rem}.
	\begin{theorem}\label{rem} Suppose that $f$ is a bounded analytic function on the right-half plane $\{z\in\mathbb{C}: \Re z>0\}$ that vanishes at a set of pairwise distinct points $d_1,d_2,\cdots,$ where
		\begin{itemize}
			\item[(i)]
			$\inf \{\vert d_n\vert\}>0$, and
			\item[(ii)]
			$\sum\limits_{n\geq 1}^{} \Re(\frac{1}{d_n})=\infty$.
		\end{itemize}
		Then $f$ vanishes identically on $\{z\in\mathbb{C}:\Re z>0\}$.
	\end{theorem}
	
	We shall often use the following classical lemma (see also \cite[Lemma 7, p.1727]{lt}).
	\begin{lemma}\label{nev}
		If a meromorphic function in a right half plane belonging to the Nevanlinna class is periodic, then it is constant.
	\end{lemma}
	
	\begin{remark}\label{periodic}
		\begin{itemize}
			\item[1)] A direct calculation shows that $\widehat{r^n}(z)=\dfrac{1}{z+n}$, for $n\in\mathbb{Z}$. 
			\item[2)] As for Theorem \ref{rem}, we will apply it in the following setting: let $(n_k)_k$ is an arithmetic sequence of positive integers, and suppose that for a certain radial function $\phi$, we have $\widehat{\phi}(n_k)=0$ for all $k$. By Theorem \ref{rem}, this implies that $\widehat{\phi}$ is identically zero on the right-half plane. Consequently, $\phi$ must also vanish on the right-half plane.
			\item[3)] We frequently rely on Lemma \ref{nev} in our arguments. In fact, in our proofs, we often encounter equations of the form  
			\[
			F(z + p) - F(z) = G(z + p) - G(z),
			\]  
			where $\Re(z) > 0$, \(p\) is an integer, and \(F\) and \(G\) are bounded analytic functions on the right half-plane. By applying Lemma \ref{nev}, we can deduce that \(F(z) = c + G(z)\) for some constant \(c\). 
			 
		\end{itemize}
	\end{remark}

	\section{Main theorem}
	Given a symbol $$g(re^{i\theta})=e^{i\theta}r^3+\sum_{l=1}^{\infty}a_l\bar{z}^l,\  a_l\in\mathbb{C},\  \textrm{ where } a_l\neq 0 \textrm{ for at least one }l\geq 5,$$ we aim to characterize all symbols of the form (i.e., symbols whose polar decomposition is truncated above)
	$$f(re^{i\theta})=\sum_{n=-\infty}^Ne^{in\theta}f_n(r),\ N\geq 1,$$ for which the associated Toeplitz operators $T_f$ commute with $T_g$.  It is understood here that $f_N\neq 0$.
	We recall that $T_f$ commutes with $T_g$ if and only if  $$T_fT_g(z^k)=T_gT_f(z^k)$$ for all vectors $z^k$ in the orthogonal basis of $L^2_a(\mathbb{D})$. Equivalently, using the notion of commutator of two operators $T$ and $S$, namely $[T,S]=TS-ST$, this condition can be expressed as
	\begin{equation}\label{commute}
		[T_f,T_g](z^k)=T_fT_g(z^k)-T_gT_f(z^k)=0, \textrm{ for all }k\geq 0.
	\end{equation}
	
	Our main theorem can be stated as follows.
	\begin{theorem}\label{main}
		Let $g$ be a symbol  of the form $g(re^{i\theta})=e^{i\theta}r^3+\sum_{l=1}^{\infty}a_l\bar{z}^l$, where $z=re^{i\theta}$, $a_l\in\mathbb{C}$ and  $a_l\neq 0$  for at least one $l\geq 5$. If there exists a nonzero function $f$ of the form $f(re^{i\theta})=\sum_{n=-\infty}^Ne^{in\theta}f_n(r)$,  with $N\geq 1,$ such that the commutator $[T_f, T_g]=0$, then $T_f$ is a polynomial of degree at most one in $T_g$. In other words, there exist constants $C_1, C_0$ such that $T_f=C_1T_g+C_0I$, where $I$ denotes the identity operator.
	\end{theorem}	
	It is important to highlight our motivation for choosing the term $e^{i\theta}r^3$ in the symbol $g$. In the work of Le and Akakai \cite[Theorem 3,~p.1725]{lt}, the analytic polynomial part of the symbol is replaced by a quasihomogeneous symbol. While this is not the most general choice, it is also not selected merely to simplify our arguments. Notably, a Toeplitz operator with the symbol $e^{i\theta}r^3$ always has powers, making the characterization of the commutant of $T_g$ significantly more challenging, as we will demonstrate in the final section. We are firmly convinced that our approach extends to more general quasihomogeneous symbols. However, such generalizations would lead to considerably more intricate calculations, which, while feasible, can become cumbersome and lengthy. We hope our choice convinces  readers of the significance of our result and the arguments, which differ entirely  from those used in the proof of \cite[Theorem 3,~p.1725]{lt}. In fact, the latter  relies on techniques that cannot be directly applied to our case, necessitating a completely different approach.
	\section{Calculations}
	\begin{lemma}\label{N/N-1}
		If equation (\ref{commute}) is satisfied, then $$T_{e^{iN\theta}f_N}=C_N\left(T_{e^{i\theta}r^3}\right)^N \textrm{ and }\  T_{e^{i(N-1)\theta}f_{N-1}}=C_{N-1}\left(T_{e^{i\theta}r^3}\right)^{N-1},$$ for some constants $C_N$ and $C_{N-1}$.
	\end{lemma}
	\begin{proof}
		In equation (\ref{commute}), the term $z^{k+N+1}$ comes only from
		\begin{equation*}
			T_{e^{i\theta}r^3}T_{e^{iN\theta}f_N}(z^k)=T_{e^{iN\theta}f_N}T_{e^{i\theta}r^3}(z^k), \textrm{ for all }k\geq 0.
		\end{equation*}
		Thus $T_{e^{iN\theta}f_N}$ commutes with $T_{e^{i\theta}r^3}$. Hence 
		\begin{eqnarray*}
			T_{e^{iN\theta}f_N}(z^k)&=&2(k+N+1)\widehat{f}_N(2k+N+2)z^{k+N}\\&=&C_N(T_{e^{i\theta}r^3})^N(z^k)\\&=&C_N\prod_{j=0}^{N-1} \frac{2(k+j+2)}{(2k+2j+6)}z^{k+N}\\&=&C_N\frac{2k+4}{2k+2N+4}z^{k+N},
		\end{eqnarray*}
		for some constant $C_N$. Here the first equality results from Lemma \ref{mellin}, the second from \cite[Proposition 7, p.1469]{l}, and the third from \cite[Lemma 3, p.1467]{l}.
		
		By similar argument we prove that $T_{e^{i(N-1)\theta}f_{N-1}}=C_{N-1}(T_{e^{i\theta}r^3})^{N-1}$, for some constant $C_{N-1}$.
	\end{proof}
	
	\begin{lemma}\label{N-2}
		If equation (\ref{commute}) is satisfied, then
		\begin{eqnarray*}\widehat{f_{N-2}}(z+N) &=& B_{N,N-2}\frac{z+4}{(z+2N-2)(z+2N)}\\&+&a_1C_N\frac{z+4}{(z+2N-2)(z+2N)}\sum_{i=0}^{N-1}\frac{z+2i}{z+2i+4},\end{eqnarray*} for some constant $B_{N,N-2}$.
	\end{lemma}
	\begin{proof}
		In equation (\ref{commute}), the term $z^{k+N-1}$ comes from 
		\begin{equation}\label{com1}
			\Big[T_{e^{i\theta}r^3},T_{e^{i(N-2)\theta}f_{N-2}}\Big](z^k)=\Big[T_{e^{iN\theta}f_{N}},T_{a_1\overline{z}}\Big](z^k), \textrm{ for } k\geq 0.
		\end{equation}
		It is understood here that $a_1\neq 0$. Thus, using Lemma \ref{mellin}, we obtain
		\begin{align*}
			&(2k+2N-2)\widehat{f_{N-2}}(2k+N)\frac{2k+2N}{2k+2N+2} -\frac{2k+4}{2k+6}(2k+2N)\widehat{f_{N-2}}(2k+N+2) \\
			&=a_1\frac{2k}{2k+2}(2k+2N)\widehat{f_{N}}(2k+N) -a_1\frac{2k+2N}{2k+2N+2}(2k+2N+2)\widehat{f_{N}}(2k+N+2)\\
			&=a_1C_N\frac{2k(2k+2)}{(2k+2)(2k+2N+2)} -a_1C_N\frac{(2k+2N)(2k+4)}{(2k+2N+2)(2k+2N+4)}.
		\end{align*}
		We complexify the equation above by taking $z=2k$, we use Remark \ref{periodic}, and we multiply both sides by $\frac{z+2N+2}{z+4}$ to obtain
		\begin{align*}
			&(z+2N-2)\widehat{f_{N-2}}(z+N)\frac{z+2N}{z+4} -\frac{z+2N+2}{z+6}(z+2N)\widehat{f_{N-2}}(z+N+2) \\
			&=a_1C_N\frac{z}{z+4} -a_1C_N\frac{z+2N}{z+2N+4} \\
			&=a_1C_N\sum_{i=0}^{N-1}\frac{z+2i}{z+2i+4} -a_1C_N\sum_{i=1}^{N}\frac{z+2i}{z+2i+4}.
		\end{align*}
		Using Remark \ref{periodic} again, we conclude that there exists a constant $B_{N,N-2}$ such that
		\begin{eqnarray*}
			\widehat{f_{N-2}}(z+N) &=& B_{N,N-2}\frac{z+4}{(z+2N-2)(z+2N)}\\&+&a_1C_N\frac{z+4}{(z+2N-2)(z+2N)}\sum_{i=0}^{N-1}\frac{z+2i}{z+2i+4}.
		\end{eqnarray*}
	\end{proof}
	Similarly, we state the following lemma, omitting the proof.
	\begin{lemma}\label{calculation}
		\begin{itemize}
			\item[i)] By considering the term $z^{k+N-2}$ that comes from
			\begin{equation}\label{f_{N-3}}
				\Big[T_{e^{i\theta}r^3},T_{e^{i(N-3)\theta}f_{N-3}}\Big](z^k)=\Big[T_{e^{iN\theta}f_{N}},T_{a_2\overline{z}^2}\Big](z^k)+\Big[T_{e^{i(N-1)\theta}f_{N-1}},T_{a_1\bar{z}}\Big](z^k)
			\end{equation}
			and assuming that $a_2, a_1$ are both nonzero, we obtain that		
			\begin{eqnarray*}
				\widehat{f_{N-3}}(z+N-1)&=&B_{N,N-3}\frac{z+4}{(z+2N-2)(z+2N-4)}\\ &+&\frac{a_2C_N(z+4)}{(z+2N-2)(z+2N-4)}\sum_{i=0}^{N-1}\frac{(z+2i-2)(z+2i)}{(z+2i+2)(z+2i+4)} \\
				&+&a_1C_{N-1}\frac{z+4}{(z+2N-2)(z+2N-4)}\sum_{i=0}^{N-2}\frac{(z+2i)}{(z+2i+4)}, 
			\end{eqnarray*}\textrm{ for some constant } $B_{N,N-3}$.
			\item[ii)] By considering the term $z^{k+N-3}$ that comes from
			\begin{eqnarray}\label{f_{N-4}}\Big[T_{e^{i\theta}r^3},T_{e^{i(N-4)\theta}f_{N-4}}\Big](z^k)&=&\Big[T_{e^{iN\theta}f_{N}},T_{a_3\overline{z}^3}\Big](z^k)+\Big[T_{e^{i(N-1)\theta}f_{N-1}},T_{a_2\bar{z}^2}\Big](z^k)\nonumber\\&+&\Big[T_{e^{i(N-2)\theta}f_{N-2}},T_{a_1\bar{z}}\Big](z^k)\end{eqnarray}
			and assuming that none of the $a_3, a_2, a_1$ is zero, we obtain that	
			\begin{eqnarray*}
				\widehat{f_{N-4}}(z+N-2)&=& B_{N,N-4}\frac{(z+4)}{(z+2N-6)(z+2N-4)} \\
				&+&\frac{a_3C_N(z+4)}{(z+2N-6)(z+2N-4)}\sum_{i=0}^{N-1}\frac{(z+2i-2)(z+2i-4)}{(z+2i+2)(z+2i+4)} \\
				&+&\frac{a_2C_{N-1}(z+4)}{(z+2N-6)(z+2N-4)}\sum_{i=0}^{N-2}\frac{(z+2i-2)(z+2i)}{(z+2i+2)(z+2i+4)} \\
				&+&\frac{a_1B_{N,N-2}(z+4)}{(z+2N-6)(z+2N-4)}\sum_{i=0}^{N-3}\frac{(z+2i)}{(z+2i+4)} \\
				&+&\frac{a_1a_1C_N(z+4)}{(z+2N-6)(z+2N-4)}\\&\times&\sum_{j=0}^{N-2}\sum_{i=0}^{N-2-j}\frac{(z+2i)(z+2i-2+2j)}{(z+2i+4)(z+2i+2+2j)}, 
			\end{eqnarray*}
			for some constant $B_{N,N-4}$.
		\end{itemize}
	\end{lemma}
	\begin{remark}\label{a_i} We draw the reader's attention to the fact that if  $a_i$ is zero, then the corresponding symbol $a_i\bar{z}^i$ is zero, and consequently, any commutator involving this term is also zero. Thus, by assuming that these $a_i$'s are nonzero, we are considering the most general case without simplifying the problem.
	\end{remark}
	
	\section{Proof of the main theorem}
	Lemmas \ref{N/N-1}, \ref{N-2}, and \ref{calculation} will guide us through the steps toward proving the main result. In the following three propositions, we will show that if equation (\ref{commute}) hold, then the value $N$ in the symbol $f$ must equal  the quasihomogeneous degree of  $e^{i\theta}r^3$ in the symbol $g$, which is $1$.
	\begin{proposition}\label{N_leq4}
		If equation (\ref{commute}) is satisfied, then $N\leq 4$.
	\end{proposition}
	
	\begin{proof}
		Suppose $N\geq 5$. By applying partial fraction decomposition to the term on the right-hand side of the expression for $\widehat{f_{N-3}}$ in Part "i)" of Lemma \ref{calculation} and using Remark \ref{periodic}, we see that the radial function $f_{N-3}$ is a polynomial in $r$ and contains the term
		\begin{equation}\label{C_N}
			a_2C_N\frac{8}{(2N-4)(2N-6)}r^{-N+3}.
		\end{equation}
		Since $N\geq 5$, it follows that $-N+3\leq -2$. Consequently, the term in (\ref{C_N}) belongs to $L^1([0,1), rdr)$ only if $a_2C_N=0$. Given that we assume $a_2\neq 0$, we must have $C_N=0$, which contradicts our assumption that $f_N\neq 0$. Therefore, we conclude that $N<5$, meaning $N\leq 4$.
	\end{proof}
	
	\begin{proposition}\label{N_leq3}
		If equation (\ref{commute}) is satisfied, then $N\leq 3 $.
	\end{proposition}
	\begin{proof} Suppose $N=4$. By taking $k=0$ in equation (\ref{f_{N-4}}), we obtain 
		\begin{equation*}
			2\widehat{f_{0}}(2) -4\widehat{f_{0}}(4) = -a_3C_4\frac{1}{5}  -a_2C_{3}\frac{3}{10} -a_1B_{4,2}\frac{1}{2}-a_1a_1C_4\frac{1}{2}\sum_{i=1}^{3}\frac{2i}{2i+4}.
		\end{equation*}
		On the other hand, from Part "ii)" of Lemma \ref{calculation}, we have 
		\begin{align*}
			&2\widehat{f_{0}}(2)= B_{4,0} +a_3C_4\sum_{i=0}^{3}\frac{(2i-2)(2i-4)}{(2i+2)(2i+4)} +a_2C_{3}\sum_{i=0}^{2}\frac{2i(2i-2)}{(2i+2)(2i+4)} \\
			&+a_1B_{4,2}\sum_{i=0}^{1}\frac{2i}{2i+4} +a_1a_1C_4\sum_{j=0}^{2}\sum_{i=0}^{2-j}\frac{2i(2i-2+2j)}{(2i+4)(2i+2+2j)}.
		\end{align*}
		and
		\begin{align*}
			&4\widehat{f_{0}}(4)= B_{4,0} +a_3C_4\sum_{i=1}^{4}\frac{(2i-2)(2i-4)}{(2i+2)(2i+4)} +a_2C_{3}\sum_{i=1}^{3}\frac{2i(2i-2)}{(2i+2)(2i+4)} \\
			&+a_1B_{4,2}\sum_{i=1}^{2}\frac{2i}{2i+4} +a_1a_1C_4\sum_{j=0}^{2}\sum_{i=1}^{3-j}\frac{2i(2i-2+2j)}{(2i+4)(2i+2+2j)}.
		\end{align*}
		Thus, the left-hand simplifies to
		\begin{equation*}
			2\widehat{f_{0}}(2)-4\widehat{f_{0}}(4) = a_3C_4-a_3C_4\frac{1}{5}  
			-a_2C_{3}\frac{3}{10} -a_1B_{4,2}\frac{1}{2} -a_1a_1C_4\frac{1}{2}\sum_{j=0}^{2}\frac{6-2j}{10-2j}.
		\end{equation*}
		Equating both sides, we conclude that  $a_3C_4=0$. Since we are assuming $a_3\neq 0$, it follows that $C_4=0$. Hence, $N\leq 3$.
	\end{proof}
	
	\begin{proposition}\label{N_leq2}
		If equation (\ref{commute}) is satisfied, then $N\leq 2 $.
	\end{proposition}
	\begin{proof}
	 Suppose $N=3$. By considering the term $z^{k-1}$ that arises from
			\begin{eqnarray}\label{f_{N-5}}\Big[T_{e^{i\theta}r^3},T_{e^{-2i\theta}f_{-2}}\Big](z^k)&=&\Big[T_{e^{3i\theta}f_{3}},T_{a_4\overline{z}^4}\Big](z^k)+\Big[T_{e^{2i\theta}f_{2}},T_{a_3\bar{z}^3}\Big](z^k)\nonumber\\&+&\Big[T_{e^{i\theta}f_{1}},T_{a_2\bar{z}^2}\Big](z^k)\nonumber\\&+&\Big[T_{f_{0}},T_{a_1\bar{z}}\Big](z^k),
			\end{eqnarray}
			and using Remark \ref{periodic}, Lemmas \ref{N/N-1}, \ref{N-2} and \ref{calculation}, while assuming that none of the $a_4, a_3, a_2, a_1$ is zero, we obtain 	
			\begin{align}\label{fminus2}
				&\widehat{f_{-2}}(z) = B_{3,-2}\frac{(z+4)}{(z-2)(z)} +a_4C_3\frac{(z+4)}{(z-2)(z)}\sum_{i=0}^{2}\frac{(z+2i-6)}{(z+2i+2)}\frac{(z+2i-4)}{(z+2i+4)}\nonumber \\
				&+a_3C_{2}\frac{(z+4)}{(z-2)(z)}\sum_{i=0}^{1}\frac{(z+2i-4)}{(z+2i+2)}\frac{(z+2i-2)}{(z+2i+4)} +\frac{a_2B_{3,1}}{(z+2)} +\frac{a_1a_1C_{2}}{(z+2)}\nonumber \\
				&+\frac{a_1a_2C_{3}(z+4)}{(z-2)(z)}\sum_{j=0}^{1}\sum_{i=0}^{1-j}\frac{(z+2i)(z+2i-2+2j)(z+2i-4+2j)}{(z+2i+4)(z+2i+2+2j)(z+2i+2j)}\nonumber \\
				&+\frac{a_1a_2C_{3}(z+4)}{(z-2)(z)}\sum_{j=0}^{1}\sum_{i=0}^{1-j}\frac{(z+2i)(z+2i-2)(z+2i-4+2j)}{(z+2i+4)(z+2i+2)(z+2i+2j)},
			\end{align} for some constant $B_{3,-2}$. Applying partial fraction decomposition to the term on the right-hand side of equation (\ref{fminus2}) and using Remark \ref{periodic}, we observe that the radial function $f_{-2}$ contains the term $(3B_{3,-2}+a_4C_3)r^{-2}$. This term belongs to $L^1([0,1),rdr)$ only if $3B_{3,-2}=-a_4C_3$. Now, by taking $k=1$ in equation (\ref{f_{N-5}}), we obtain
		\begin{equation*}
			\widehat{f_{-2}}(4)= a_4C_3\frac{1}{15} +a_3C_{2}\frac{1}{10} +a_2B_{3,1}\frac{1}{6}+\frac{1}{6}a_1a_1C_2+(\frac{1}{5}+\frac{1}{18}+\frac{1}{12})a_1a_2C_3.
		\end{equation*}
		On the other hand, equation (\ref{fminus2}) implies
		\begin{equation*}
			\widehat{f_{-2}}(4) = B_{3,-2}+a_4C_3\frac{1}{15} +a_3C_{2}\frac{1}{10} +a_2B_{3,1}\frac{1}{6}+\frac{1}{6}a_1a_1C_2+(\frac{1}{5}+\frac{1}{18}+\frac{1}{12})a_1a_2C_3.
		\end{equation*}
		Thus, we must have $B_{3,-2}=0$. Since $3B_{3,-2}=-a_4C_3$ and  we are assuming $a_4\neq 0$, it follows that $C_3=0$. Hence, $N\leq 2$.
	\end{proof}
	
	\begin{proposition}\label{N_leq1}	If equation (\ref{commute}) is satisfied, then $N\leq 1 $. Moreover, under the assumption that $N\geq 1$, we conclude that $N=1$.
	\end{proposition}
	
	\begin{proof}
		Suppose $N=2$. In this case, Part "i)" of Lemma \ref{calculation} implies that
		\begin{eqnarray}\label{fminus1}
			\widehat{f_{N-3}}(z+N-1)&=&\widehat{f_{-1}}(z+1)\nonumber\\ &=&B_{2,-1}\frac{z+4}{z(z+2)} +\frac{z-2}{(z+2)(z+2)}+\frac{a_2C_2}{z+6} +\frac{a_1C_{1}}{z+2}.
		\end{eqnarray}
		Now, by taking $k=0$ in equation (\ref{f_{N-3}}), we obtain
		\begin{equation*}
			\widehat{f_{-1}}(3) = a_2C_2\frac{1}{8}+ a_1C_{1}\frac{1}{4}.
		\end{equation*}
		On the other hand, equation (\ref{fminus1}) yields
		\begin{equation*}
			\widehat{f_{-1}}(3) =B_{2,-1}\frac{3}{4} +a_2C_2\frac{1}{8} +a_1C_{1}\frac{1}{4}.
		\end{equation*}
		Thus, we must have $B_{2,-1}=0$. Next, we analyze the term $z^{k-1}$, which arises from
			\begin{eqnarray}\label{f_{-2}}\Big[T_{e^{i\theta}r^3},T_{e^{-2i\theta}f_{-2}}\Big](z^k)&=&\Big[T_{e^{2i\theta}f_{2}},T_{a_3\overline{z}^3}\Big](z^k)+\Big[T_{e^{i\theta}f_{1}},T_{a_2\bar{z}^2}\Big](z^k)\nonumber\\&+&\Big[T_{f_{0}},T_{a_1\bar{z}}\Big](z^k).\end{eqnarray}
			Assuming that none of the $a_3, a_2, a_1$ is zero, Remark \ref{periodic} and Lemmas \ref{N/N-1} and  \ref{N-2} imply
			\begin{eqnarray*}
				\widehat{f_{-2}}(z)&=& B_{2,-2}\frac{(z+4)}{(z-2)(z)} + \frac{a_3C_2(z+4)}{(z-2)(z)}\sum_{i=0}^{1}\frac{(z+2i-2)(z+2i-4)}{(z+2i+2)(z+2i+4)} \\
				&+&\frac{a_2C_{1}}{(z+2)} + \frac{a_1a_1C_2}{(z+2)}, 
			\end{eqnarray*}
			for some constant $B_{2,-2}$.
		By applying partial fraction decomposition to the term on the right-hand side of the equation above and using Remark \ref{periodic}, we see that the radial function $f_{-2}$ contains the term $3B_{2,-2}r^{-2}$. Therefore, for $f_{-2}$ to be in $L^1([0,1),rdr)$, we must have $B_{2,-2}=0$. Finally, we consider the term $z^{k-3}$, which comes from
			\begin{eqnarray*}\Big[T_{e^{i\theta}r^3},T_{e^{-4i\theta}f_{-4}}\Big](z^k)&=&\Big[T_{e^{2i\theta}f_{2}},T_{a_5\overline{z}^5}\Big](z^k)+\Big[T_{e^{i\theta}f_{1}},T_{a_4\bar{z}^4}\Big](z^k)\\&+&\Big[T_{f_{0}},T_{a_3\bar{z}^3}\Big](z^k)+\Big[T_{e^{-i\theta}f_{-1}},T_{a_2\bar{z}^2}\Big](z^k)\\&+&\Big[T_{e^{-2i\theta}f_{-2}},T_{a_1\bar{z}}\Big](z^k).
			\end{eqnarray*} 
	 Assuming none of the $a_5, a_4, a_3, a_2, a_1$ is zero, Remark \ref{periodic} and Lemmas \ref{N/N-1}, \ref{N-2} and \ref{calculation} imply
			\begin{eqnarray*}
				\widehat{f_{-4}}(z-2) &=& B_{2,-4}\frac{(z+4)}{(z-6)(z-4)} \\
				&+&a_5C_2\frac{(z+4)}{(z-6)(z-4)}\sum_{i=0}^{1}\frac{(z+2i-8)}{(z+2i+2)}\frac{(z+2i-6)}{(z+2i+4)} \\
				&+&\frac{a_4C_{1}}{(z+2)}  +2\frac{a_1a_3C_2}{(z+2)} + \frac{a_2a_2C_2}{(z+2)},
			\end{eqnarray*}
			for some  constant  $B_{2,-4}$. By applying partial fraction decomposition to the term on the right-hand side of the equation above and using Remark \ref{periodic}, we clearly see that the radial  function $f_{-4}$ contains the terms $5B_{2,-4}r^{-4}$ and $(-4B_{2,-4}-\frac{2}{3}a_5C_2)r^{-2}$. Therefore, for $f_{-4}$ to be in  $L^1([0,1),rdr)$, we must have  $B_{2,-4}=0$ and $-4B_{2,-4}-\frac{2}{3}a_5C_2=0$. Consequently,  $a_5C_2=0$. Since we are assuming $a_5\neq 0$, it follows that $C_2=0$. Hence, $N\leq 1$.
	\end{proof}	
The following lemma will be used to explicitly determine the radial functions $f_{n}$ in the symbol $f$, for $n\leq -1$. 
\begin{lemma} \label{lem1}
Let $l\in\mathbb{N}$. Assume  
\begin{equation}\label{eq_lem1_1}
\Big[T_{e^{i\theta}r^3},T_{e^{-il\theta}f_{-l}}\Big](z^k)=\Big[cT_{e^{i\theta}r^3},T_{a_l\overline{z}^l}\Big](z^k),
\end{equation}
for some constant $c$. Then $T_{e^{-il\theta}f_{-l}}=cT_{a_l\overline{z}^l}$.

\end{lemma}
\begin{proof}
If equation (\ref{eq_lem1_1}) holds, then for all $k\geq l$, we have
\begin{align*}
&(2k-2l+2)\widehat{f_{-l}}(2k-l+2)\frac{(2k-2l+4)}{(2k-2l+6)} -\frac{(2k+4)}{(2k+6)}(2k-2l+4)\widehat{f_{-l}}(2k-l+4) \\
&= ca_l\frac{(2k-2l+2)}{(2k+2)}\frac{(2k-2l+4)}{(2k-2l+6)} -ca_l\frac{(2k+4)}{(2k+6)}\frac{(2k-2l+4)}{(2k+4)}.
\end{align*}
We complexify the equation above by taking $z=2k$. Using Remark \ref{periodic}, we multiply both sides by $\frac{(z-2l+6)}{(z+4)}$ to obtain
\begin{align*}
&(z-2+2)\widehat{f_{-l}}(z-l+2)\frac{(z-2l+4)}{(z+4)} -\frac{(z-2l+6)}{(z+6)}(z-2l+4)\widehat{f_{-l}}(z-l+4) \\
&= ca_l\frac{(z-2l+2)(z-2l+4)}{(z+2)(z+4)} -ca_l\frac{(z-2l+4)(z-2l+6)}{(z+4)(z+6)}.
\end{align*}
Using Remark \ref{periodic} again, we find that there exists a constant $A_l$ such that
 \begin{equation}\label{fminusl}
\widehat{f_{-l}}(z-l+2) =A_l\frac{(z+4)}{(z-2l+4)(z-2l+2)} + ca_l\frac{1}{(z+2)} .
\end{equation}
We now consider the following cases:
\begin{itemize}
\item[ \textbf{Case \( l\geq 2. \)}] By applying partial fraction decomposition to the term on the right-hand side of equation (\ref{fminusl}) and using Remark \ref{periodic}, we see that the radial function $f_{-l}$  contains the term $A_l(l+1)r^{-l}$. Since $-l\leq -2$,  the function $f_{-l}$ belongs to $L^1([0,1),rdr)$ only if $A_l=0$. Hence, $f_{-l}(r)=ca_lr^{l}$. 

\item[\textbf{Case \( l=1. \)}] By taking $k=0$ in the equation (\ref{eq_lem1_1}), we obtain
\begin{equation*}
 \widehat{f_{1}}(3) = ca_1\frac{1}{4}.
\end{equation*}
On the other hand, equation (\ref{fminusl}) implies $$\widehat{f_{1}}(3) =A_1\frac{3}{4} + ca_1\frac{1}{4}.$$ Therefore, we conclude that $A_1=0$, which implies that $f_{-1}(r)=ca_1r$.

\end{itemize}
\end{proof}
Finally, we reach the last step of the proof, which is to reconstruct the symbol $f$.
\begin{proposition}
If equation (\ref{commute}) holds, then $f(re^{i\theta})=C_1g(re^{i\theta})+C_0$, for some constants $C_1$ and $C_0$.
\end{proposition}
\begin{proof}
If the equation (\ref{commute}) holds, then Proposition \ref{N_leq1} implies that $$f(re^{i\theta})=\sum_{n=-\infty}^{1}e^{in\theta}f_n.$$\\
From Lemma \ref{N/N-1}, we have that $f_1(r)=C_1r^3$ and $f_0=C_0$.
Next, for all integers $l\geq 1$, the term $z^{k-l}$ comes only from
\begin{equation*}
\Big[T_{e^{i\theta}r^3},T_{e^{-il\theta}f_{-l}}\Big](z^k)=\Big[C_1T_{e^{i\theta}r^3},T_{a_l\overline{z}^l}\Big](z^k).
\end{equation*}
Thus, Lemma \ref{lem1} implies that $f_{-l}=C_1r^l$.
Finally, by combining those  results, we obtain
		\begin{eqnarray*}f(re^{i\theta})&=&e^{i\theta}f_1(r)+f_0(r)+\sum_{n=-\infty}^{-1}e^{in\theta}f_n(r)\\
			&=& C_1 e^{i\theta}r^3+C_0+\sum_{l=1}^{\infty}C_1a_le^{-il\theta}r^l\\
			&=& C_1g(re^{i\theta})+C_0.
			\end{eqnarray*}
			Using the fact that Toeplitz operators are linear with respect to their symbols, we conclude that
			$$T_f=C_1T_g+C_0I.$$
\end{proof}		
\section{Special case} In this section, we consider the case where the infinite anti-analytic power series in the symbol $g$ from Theorem \ref{main} is replaced by a polynomial in $\bar{z}$.
		\begin{lemma}
			Let $g(re^{i\theta})=e^{i\theta}r^3+\sum\limits_{l=1}^{m}a_l e^{-il\theta}r^l$. The 
			the product $T_{g}^2$ is Toeplitz operator if and only if $m\leq 4$.
		\end{lemma}
		\begin{proof}
			We have
			\begin{align*}
				T_{g}^2 &=T_{e^{i\theta}r^3+\sum\limits_{l=1}^{m}a_l e^{-il\theta}r^l}T_{e^{i\theta}r^3+\sum\limits_{l=1}^{m}a_l e^{-il\theta}r^l} \\
				&= T_{e^{i\theta}r^3}^2+\sum\limits_{l=1}^{m}a_lT_{e^{i\theta}r^3}T_{\overline{z}^l}+\sum\limits_{l=1}^{m}a_lT_{\overline{z}^l}T_{e^{i\theta}r^3}+\sum\limits_{l=1}^{m}a_lT_{\overline{z}^l}T_{\overline{z}^l}  \\
				&= T_{e^{i\theta}r^3}^2+\sum\limits_{l=1}^{m}a_lT_{e^{i\theta}r^3}T_{\overline{z}^l}+\sum\limits_{l=1}^{m}a_lT_{\overline{z}^l}T_{e^{i\theta}r^3}+\sum\limits_{l=1}^{m}a_lT_{\overline{z}^{2l}}.
			\end{align*}
			It is well-known (see \cite[Example 3.8,~p.1057]{bbl}) that $T_{{e^{i\theta}r^3}}^n$ is always a Toeplitz operator for all $n\geq 1$. Because $\bar{z}^l$ is an anti-analytic symbol,  $T_{\bar{z}^l}T_{e^{i\theta}r^3}=T_{e^{i\theta-l}r^{3+l}}$ for all $l\geq 1$. Thus, we focus  on the product \( T_{e^{i\theta}r^3}T_{\overline{z}^l} \). For all \( k \geq l \), we have 	
			\begin{align*}
				T_{e^{i\theta}r^3}T_{\overline{z}^l}(z^k) &= \frac{(2k-2l+2)}{(2k+2l+2)}\frac{(2k-2l+4)}{(2k-2l+6)} z^{k+1-l} \\
				&= (2k-2l+4)\widehat{h}(2k-l+3)z^{k+1-l},
			\end{align*}
			where the radial function $h$ is defined as follows:
			\begin{itemize}
				\item[$\bullet$] $h(r)=\frac{l}{l-1}r^{3l-1}+\frac{1}{1-l}r^{3-l}$ if $l\geq 2$. In this case, $h$ belongs to $L^1([0,1),rdr)$ if and only if $3-l+1\geq 0$, which simplify to $l\leq 4$.
				\item[$\bullet$] $h(r)=r^{2}+4r^{2}\log r$ if $l=1$. Here, $h$ is in $L^1([0,1), rdr)$.
			\end{itemize}
	Therefore, for $m\leq 4$, $T^2_g$ remains a Toeplitz operators, completing the proof.
					\end{proof}
The above lemma allows us to replace the anti-analytic part of the symbol $g$ in Theorem \ref{main} by a polynomial in $\bar{z}$ of degree at most $4$, yielding the following result. We provide a brief proof, as the calculations are nearly identical to those in the proof of our main theorem in the previous sections.
\begin{theorem}
		Let $g$ be a symbol  of the form $g(re^{i\theta})=e^{i\theta}r^3+\sum_{l=1}^{m}a_l\bar{z}^l$, where $z=re^{i\theta}$, $a_l\in\mathbb{C}$ and $m\leq 4$. If there exists a nonzero function $f$ of the form $f(re^{i\theta})=\sum_{n=-\infty}^Ne^{in\theta}f_n(r)$,  with $N\geq 1,$ such that the commutator $[T_f, T_g]=0$, then $T_f=P(T_g)$, where $P$ is a polynomial of degree at most two.
	\end{theorem}				
\begin{proof}
	Following the same argument as in the previous sections, we deduce that if $[T_f,T_g]=0$, then $N\leq 2$. Finally, using the same techniques of the proof of \cite[Theorem 2,~p. 886]{lry}, we conclude that $T_f$ is a polynomial of degree at most two in $T_g$.
	\end{proof}

\end{document}